\documentclass[12pt]{amsart}
\usepackage{amstext,amsfonts,amsmath,amssymb,amsthm}

\newfont{\gl}{eufm10 scaled \magstep1} 

\newcommand{\dd}{\hbox{d}}

\def\Ker{{\rm Ker}}

\numberwithin{equation}{section}

\newtheorem{theorem}{Theorem}[section]

\newtheorem{proposition}{Proposition}[section]

\newtheorem{lemma}{Lemma}[section]

\theoremstyle{remark}
\newtheorem{example}{Example}[section]

\theoremstyle{remark}
\newtheorem{remark}{Remark}[section]

\theoremstyle{definition}
\newtheorem{definition}{Definition}[section]

\newcommand {\un}{\underline}

\begin{document}

\title{{Deformation of Dirac structures\\ along isotropic
subbundles}}

\author{Iv\'an Calvo} 
\address{Laboratorio Nacional de Fusi\'on\\
Asociaci\'on EURATOM-CIEMAT\\
E-28040 Madrid (Spain)}
\email{ivan.calvo@ciemat.es} 

\author{Fernando Falceto}
\address{Departamento de F\'{\i}sica Te\'orica and
Instituto de Biocomputaci\'on y F\'{\i}sica de Sistemas Complejos,
Universidad de Zaragoza,
E-50009 Zaragoza (Spain)}
\email{falceto@unizar.es}

\author{Marco Zambon}
\address{Universidade do Porto, Departamentos de Matematica Pura, Rua do Campo Alegre 687, 4169-007 Porto (Portugal)}
\email{mzambon@fc.up.pt}

\date{\today}

\begin{abstract}
\noindent
Given a Dirac subbundle and an isotropic subbundle of a Courant algebroid, we provide a
canonical method to obtain a new Dirac subbundle. When the original
Dirac subbundle is involutive (i.e., a Dirac structure) this construction has
interesting applications, for instance  to Dirac's theory of constraints and to the Marsden-Ratiu reduction in Poisson geometry.
\end{abstract}
\maketitle

\section{Introduction} 
The concept of {\it Dirac
structure}  generalizes  Poisson and presymplectic structures
by embedding them in the framework of the geometry of $TM\oplus
T^*M$ or, more generally, the geometry of a \emph{Courant algebroid}. Dirac structures were introduced in a remarkable paper by
T. Courant \cite{Cou}.  Therein, they are related to the
Marsden-Weinstein reduction \cite{MarWei} and to the Dirac bracket
\cite{Dir} on a submanifold of a Poisson manifold. More recently, Dirac
structures have been considered in connection to the reduction of
implicit Hamiltonian systems (see \cite{BlaSch},\cite{BlaRat}).  This
simple but powerful structure allows to deal with mechanical
situations in which we have both gauge symmetries and Casimir
functions.

We present a construction
 which takes an isotropic subbundle $S$ and a Dirac subbundle $D$ of an exact Courant algebroid, and  produces a new Dirac subbundle $D^S$ (Def. \ref{stretch}). 
This construction, which we refer to as \emph{stretching}, was introduced by the first two authors in \cite{CaFa}.
 When both $S$ and $D$ are involutive, we find conditions ensuring that $D^S$ is also involutive, i.e. a Dirac structure (Thm. \ref{connection}).  

We further show that three  prominent classes of Dirac structures are indeed stretched Dirac structures: the Dirac brackets that appeared in Dirac's theory of constraints, the Dirac structures underlying the Marsden-Ratiu quotients in Poisson geometry \cite{MarRat}, and coupling Dirac structures on Poisson fibrations \cite{BF}. 

The paper is organized as follows. In Section 2 we review basic definitions, in Section 3 we describe our stretching construction, in Section 4 we discuss when the stretched structure is involutive, and in Section 5 we present examples and applications.

\section{Courant algebroids and Dirac structures}

\begin{definition}\label{defCA} A \textbf{Courant algebroid} \cite{LWX} over a
manifold $M$ is a vector bundle $E\rightarrow M$ equipped with an
${\mathbb R}$-bilinear bracket $[\cdot,\cdot]$ on $\Gamma(E)$, a
non-degenerate symmetric bilinear form $\langle\cdot,\cdot\rangle$ on
the fibers and a bundle map $\pi:E\rightarrow TM$ (the {\it anchor})
satisfying, for any $e_1,e_2,e_3\in\Gamma(E)$ and $f\in C^\infty(M)$:
\begin{itemize}
\item[(i)]$[e_1,[e_2,e_3]]=[[e_1,e_2],e_3]+[e_2,[e_1,e_3]]$
\item[(ii)]$\pi([e_1,e_2])=[\pi(e_1),\pi(e_2)]$
\item[(iii)]$[e_1,fe_2]=f[e_1,e_2]+(\pi(e_1)f)e_2$
\item[(iv)]$\pi(e_1)\langle e_2,e_3\rangle=\langle [e_1,e_2],e_3\rangle+\langle e_2,[e_1,e_3]\rangle$
\item[(v)]$[e,e]={\mathcal D}\langle e,e\rangle$
\end{itemize}
where ${\mathcal  D}:C^\infty(M)\rightarrow \Gamma(E)$ is defined by
${\mathcal  D}=\frac{1}{2}\pi^*\circ \dd$, using the bilinear form to identify
$E$ and its dual.
\end{definition}

We see from axiom (v) that the bracket is not skew-symmetric,
but rather satisfies
$
[e_1,e_2]=-[e_2,e_1]+2{\mathcal  D} \langle e_1,e_2 \rangle.$

A Courant algebroid is called \textbf{exact} if
\begin{equation*}
0\longrightarrow
T^*M\stackrel{\pi^*}{\longrightarrow}E\stackrel{\pi}{\longrightarrow}TM\longrightarrow
0
\end{equation*}
is an exact sequence. 
Choosing a splitting $TM \rightarrow E$ of the above sequence with isotropic image allows one to identify
the exact Courant algebroid   with
$TM\oplus T^*M$ endowed with the natural symmetric pairing
\begin{equation*}
\langle(X,\xi),(X',\xi')\rangle = \frac{1}{2}(i_{X'}\xi + i_{X}\xi')
\end{equation*}
and the Courant bracket
\begin{equation*}\label{Courantbracket}
[(X,\xi),(X',\xi')]=([X,X'], {\mathcal  L}_{X}\xi'-i_{X'}\dd\xi  +i_{X'}i_{X}H)
\end{equation*}
for some closed 3-form $H$. In fact, the Courant algebroid uniquely
determines the cohomology class of $H$, called \v{S}evera class. The anchor $\pi$ is given by
the projection onto the first component. When it is important to
stress the value of the 3-form $H$ we shall use the
notation $E_H$ for $TM\oplus T^*M$ equipped with this Courant
algebroid structure.
 
\begin{definition}A \textbf{Dirac subbundle} or \textbf{almost Dirac structure}  in an exact Courant algebroid $E$ is a subbundle $D\subset E$ which is
maximal isotropic with respect to
$\langle\cdot,\cdot\rangle$. The maximal isotropicity condition implies that $D^\perp =
D$, where $D^\perp$ stands for the orthogonal subspace of $D$. In
particular, $\mbox{rank}(D)=\mbox{dim}(M)$.

A \textbf{Dirac structure} is an involutive Dirac subbundle, i.e. a
Dirac subbundle $D$ whose sections closed under the Courant bracket.  In
this case the restriction to $D$ of the Courant bracket is
skew-symmetric and $D$ with anchor $\pi$ is a Lie algebroid.
\end{definition}

The two basic examples of Dirac structures are:

\begin{example}
For any 2-form $\omega$, the graph $L_\omega$ of
$\omega^\flat:TM\rightarrow T^*M$ is a Dirac subbundle such that
$\pi(L_\omega)=TM$. $L_\omega$ is a Dirac
structure in $E_H$ if and only if $\dd\omega =-H$. In particular,
$L_\omega$ is a Dirac structure in $E_0$ if and only if $\omega$ is
closed.
\end{example}

\begin{example}
Let $\Pi$ be a bivector field on $M$. The graph $L_\Pi$ of the map
$\Pi^\sharp:T^*M\rightarrow TM$ is always a Dirac subbundle. In this
case the natural projection from $L_\Pi$ to $T^*M$ is
one-to-one. $L_\Pi$ is a Dirac structure in $E_H$ if and only if $\Pi$
is a twisted Poisson structure. In particular, $L_\Pi$ is a Dirac
structure in $E_0$ if and only $\Pi$ is a Poisson structure.
\end{example}

%
 
  
\section{Stretched Dirac subbundles}
 
Until the end of this note we will assume the following setup:
\begin{center}
\fbox{
\parbox[c]{12.6cm}{
\begin{center}
$E$ is an exact
Courant algebroid\\
$S\subset E$ is an isotropic subbundle (i.e. $S\subset S^{\perp}$)\\
$D\subset E$ is  an Dirac subbundle.\\ 
\end{center}
}}
 \end{center}
We further assume that $D\cap
S$ (or equivalently $D\cap S^{\perp}$) has constant rank along $M$.

\begin{definition}\label{stretch}
The \textbf{stretching of $D$ along $S$} \cite{CaFa} is the Dirac subbundle
\begin{equation}\notag
D^S:= (D\cap
S^\perp) + S.
\end{equation}
\end{definition}

To justify the fact that $D^S$ is maximal
isotropic we use
\begin{eqnarray}\notag
(D^S)^\perp &=& (D^\perp + S)\cap S^\perp\cr
&=& (D\cap S^\perp)+S =D^S,
\end{eqnarray}  
where in the last line we have used that $D$ is maximal isotropic
and $S$ is a subset of $S^\perp$. It is also clear that $D^S$, as the
sum of two subbundles whose intersection  has constant rank, is a (smooth) subbundle.

$D^S$ is the Dirac subbundle
closest to $D$ among those containing $S$, as stated in the following
\begin{proposition}\cite{CaFa}\label{close}
Let $D,S$ and $D^S$ be as above and let $D'$ be a Dirac subbundle such
that $S\subset D'$. Then, $ D'\cap D\subset D^S\cap D$. In addition,
$D'\cap D =D^S\cap D$ if and only if $D'=D^S$.
\end{proposition}
\begin{proof}
From the isotropicity of $D'$ and $S\subset D'$ we deduce
 that $D'\subset S^{\perp}$. Hence,
$$D'\cap D \subset S^\perp\cap D = D^S\cap D.$$

If the equality $D'\cap D =D^S\cap D$ holds, then $D'\supset D'\cap D
= S^\perp\cap D$. Since $S\subset D'$, we find that $D^S=(D\cap
S^\perp) + S\subset D'$. But $D^S$ and $D'$ have the same dimension,
so that they are equal. 
\end{proof}

\section{Integrability}

In this section we  determine various properties of $D^S$, in particular conditions under which $D^S$ is a Dirac structure.

\begin{lemma}\label{th:invariance} Assume   that  $S$ and $D$ are closed under the Courant bracket. Then the set of $S$-invariant sections of $D^S$
$$ \{e\in D^S: [\Gamma(S),e]\subset \Gamma(S) \}$$ is closed under the Courant
bracket.
\end{lemma}
\begin{proof}
Consider two sections $e_1,e_2\in \Gamma(D^S)$ which are $S$-invariant, i.e.,   $[\Gamma(S),e_i]\in \Gamma(S)$. First, let us prove that
$[e_1,e_2]$ is an $S$-invariant section. Take $s\in\Gamma(S)$ and
write
\begin{eqnarray}\notag
[s,[e_1,e_2]]&=&[[s,e_1],e_2]+
[e_1, [s, e_2]]\end{eqnarray}
by Def. \ref{defCA} i). Now recall that $[e,s]= -[s,e]$ for $e\in \Gamma(D^S)$
and $s\in\Gamma(S)$ because $D^S=(D+S)\cap S^\perp\subset
S^\perp$. The $S$-invariance of $[e_1,e_2]$ follows immediately.

Next we show that $[e_1,e_2]\in\Gamma(D^S)$. Since we assumed that both
$D\cap S^\perp$ and $S$ are subbundles, every section $e\in
\Gamma(D^S)$ can be written as $e=v+w$ with $v\in\Gamma(D\cap
S^\perp)$ and $w\in \Gamma(S)$. Notice that if $e$ is $S$-invariant,
$v$ is also $S$-invariant because $S$ is Courant involutive. The
expression
\begin{equation}\label{eq:proofDplusS}
[e_1,e_2]=[v_1+w_1,v_2+w_2]=[v_1,v_2]+[v_1,w_2]+[w_1,v_2]+[w_1,w_2]
\end{equation}
makes clear that $[e_1,e_2]\in\Gamma(D+S)$, since
$[v_1,v_2]\in\Gamma(D)$ and the remaining terms on the right-hand side
of eq. \eqref{eq:proofDplusS} are sections of $S$. To prove
that $[e_1,e_2]\in\Gamma(S^\perp)$ notice that, for any $s\in\Gamma(S)$,
\begin{eqnarray}\notag
\langle s,[e_1,e_2]\rangle &=&
\pi(e_1)\langle s,e_2\rangle-\langle [e_1,
s],e_2\rangle=0
\end{eqnarray}
where we have used Def. \ref{defCA} iv)
as well as  the orthogonality of $s$ and $e_i,\ i=1,2$. 
\end{proof}
   
Inspired by \cite{OrtRat} we will give the following
\begin{definition}\label{canonical}
Given a Dirac subbundle $D$ and an involutive isotropic 
 subbundle $S\subset E$, 
we say that $S$ is \textbf{canonical} for $D$ if there exists a local 
$S$-invariant section\footnote{Recall that a section $e$ is $S$-invariant iff 
 $[\Gamma(S),e]\subset \Gamma(S)$.}
 of $D^S$ passing through any 
of its points. 
\end{definition}
 
\begin{proposition}
\label{imply}
Assume   that  $S$ and $D$ are closed under the Courant bracket. We have the following chain of implications:
 \begin{itemize}
\item [a)] $S$ is canonical for $D$ $\Rightarrow$
\item [b)]$D^S$ is  a Dirac structure $\Rightarrow$
\item [c)]$[\Gamma(S),\Gamma(D^S)]\subset\Gamma(D^S)$ (i.e. $S$ preserves $D^S$)
\end{itemize}
\end{proposition}
\begin{proof} a) $\Rightarrow$ b): We have to show that the Courant bracket of two sections
$v,v'\in\Gamma(D^S)$ is again a section of $D^S$. We
write $v$ and $v'$ in terms of a local basis, $\{e_i\}$, of $S$-invariant
sections (such a basis always exists due to the canonicity of $S$). From Def. \ref{defCA} iii) and Lemma \ref{th:invariance} one
immediately obtains that $[v,v']$ is a linear combination of the
$e_i$'s and hence belongs to $\Gamma(D^S)$.
 
b) $\Rightarrow$ c): holds because $S\subset D^S$.
\end{proof}

The following theorem 
gives sufficient conditions to ensure that 
$D^S$ is a Dirac structure.  
 \begin{theorem}\label{connection}
Assume   that  $S$ and $D$ are closed under the Courant bracket and additionally that $\pi(S^\perp)$ is an
integrable regular distribution. Then items a),b),c) of Prop. \ref{imply} are all equivalent. 
\end{theorem}

\begin{proof} We just need to show c) $\Rightarrow$ a) in Prop. \ref{imply}. Notice first that $\pi(S)$ is a regular
integrable distribution. Indeed
$$\Ker(\pi)\cap S = \pi^*(\pi(S^\perp) ^{\circ}),$$ and given that
$\pi(S^\perp)$ is regular and $\pi^*$ is injective for exact Courant
algebroids, it follows that $\Ker(\pi)\cap S$ is a subbundle. Now, the fact that
$S$ is also a subbundle implies the regularity of
$\pi(S)$. Integrability follows from the assumption that $\Gamma(S)$ is closed under the Courant bracket.

Take a commuting  basis of local sections of
$\pi(S)$ denoted by $\{\partial_i\}$. Fix lifts $s_i$
 of $\partial_i$ to $S$, i.e. $s_i\in\Gamma(S)$ and
$\pi(s_i)=\partial_i$. Since we are assuming c) of Prop. \ref{imply}
we can
define  a partial $\pi(S)$-connection on $D^S$  by imposing
$$\nabla_i e :=[s_i,e].$$

The involutivity of $S$ and Def. \ref{defCA} ii) imply  that 
$[\Gamma(S),\Gamma(\Ker(\pi)\cap S)]\subset\Gamma(\Ker(\pi)\cap S)$,
so we can use $\nabla$ to define a partial $\pi(S)$-connection\footnote{The connection $\widetilde\nabla$ depends on the choice of lifts $s_i\in \Gamma(S)$.}
 $\widetilde\nabla$ on
$D^S/(\Ker(\pi)\cap S)$. We now argue that $\widetilde\nabla$  is flat.

The curvature of the connection $\nabla$, with components $F_{ij}$, is given by
\begin{equation}\label{curv}
F_{ij} e =\nabla_i \nabla_j e -\nabla_j \nabla_i e=
[s_i[s_j,e]]-[s_j[s_i,e]]
= [[s_i,s_j],e],
\end{equation}
and given that $\partial_i$ and $\partial_j$ commute
and $S$ is involutive we have 
$[s_i,s_j]\in \Ker(\pi)\cap S$.

Next we want to show that 
\begin{eqnarray}\label{Kerninvar}
[\Gamma(\Ker(\pi)\cap S),\Gamma(D^S)]\subset\Gamma(\Ker(\pi)\cap S).
\end{eqnarray}
For that, take a section $s\in\Gamma(\Ker(\pi)\cap S)$ and write  
$s=\pi^*(\eta)$ with $\eta\in \Gamma(\pi(S^\perp) ^{\circ})$. Also take
arbitrary sections $e\in\Gamma(D^S)$ and $s^\perp\in\Gamma(S^\perp)$.
Now
\begin{eqnarray*}\notag
\langle [s,e],s^\perp \rangle &=& 
\langle\pi^*(\eta), [e,s^\perp]\rangle- 
{\pi(e)}  \langle s,s^\perp\rangle\cr
&=& i_{\pi([e,s^\perp])}\eta\cr &=&i_{[\pi(e),\pi(s^\perp)]}\eta\\
&=&0,
\end{eqnarray*}
where in the first equality we used Def. \ref{defCA} iv) and in
last equality we  used that $D^S\subset S^\perp$ and
$\pi(S^\perp)$ is integrable.


Eq. \eqref{curv} and eq. \eqref{Kerninvar} together imply that  $\widetilde\nabla$ is a flat
connection.   Hence through any point of $D^S/(\Ker(\pi)\cap S)$ passes
a  local horizontal sections for $\widetilde\nabla$, and any lift of it to a   
 section $e_h$ of $D^S$
satisfies $\nabla_i e_h\in\Gamma(\Ker(\pi)\cap S)$ for all $i$. But the sections
$\{s_i\}$ used to build the connection $\nabla$, together with
$\Gamma(\Ker(\pi)\cap S)$, span $\Gamma(S)$. Hence  using Def. \ref{defCA} iii) and  eq. \eqref{Kerninvar} 
we get
$$[\Gamma(S),e_h]\subset\Gamma(S),$$
completing the proof.
\end{proof}

\begin{remark}\label{preserved}
With the hypotheses of Theorem \ref{connection}, if   $[\Gamma(S),\Gamma(D)]\subset\Gamma(D)$ then $S$ is canonical for $D$ (because condition c) in Prop. \ref{imply} is satisfied). The converse is not true, see e.g. \cite{OrtRat} for a counterexample
in the context of Poisson manifolds. \end{remark}

\section{Examples and applications}

In this section we work only with the exact Courant algebroid  $E_0$.
The first two examples show that two well-known constructions in Poisson geometry, once phrased in tensorial terms, correspond to the stretching of Poisson structures.

\subsection {Dirac brackets} We give a description of the classical Dirac bracket in tensorial terms, i.e. in terms of Dirac structures, thereby giving a clear geometric
 interpretation to the   Dirac bracket. Further we  present a natural generalization.

 
We recall first the construction of the Dirac bracket on a Poisson manifold
$(M,\Pi)$.  
\begin{definition}\label{Dirbr}
Given a regular foliation $R$ on an open set
$U\subset M$ whose leaves $N$ are
the level
sets of  second class constraints $\varphi^1,\dots,
\varphi^m$ (i.e. independent functions for which the matrix
$C^{ab}:=\{\varphi^a,\varphi^b\}_{\Pi}$ is invertible, with inverse
$C_{ab}$), the {\bf Dirac bracket} is defined as
\begin{equation} \label{dir}\{f,g\}_{Dirac}:=\{f,g\}_{\Pi}-\{f, \varphi^a\}_{\Pi}C_{ab}\{ \varphi^b,g\}_{\Pi}.
\end{equation}
We denote by $\Pi_{Dirac}$ the bivector field corresponding to the bracket $\{\cdot,\cdot\}_{Dirac}$.
\end{definition}

\begin{lemma}\label{PAM}
i) The level sets $N$ of $\varphi$ are \emph{cosymplectic} submanifolds of $(M,\Pi)$ 
and therefore have a Poisson structure induced by $\Pi$. 

ii)   $(M,\Pi_{Dirac})$ is obtained putting together the level sets $N$ of $\varphi$, endowed with the Poisson structure induced by $\Pi$. In particular the Dirac bracket \eqref{dir} depends only on the level sets of the constraints (and not on the constraints themselves).
\end{lemma}
\begin{remark}
1) Here and in the following we use repeatedly the following fact: a Poisson (Dirac) manifold is 
determined by its foliation into symplectic (presymplectic) leaves.\\
2) Lemma \ref{PAM} ii) recovers the fact that \eqref{dir} is a Poisson bracket. 
\end{remark}
\begin{proof}
i)  Since the $\varphi^i$ are second class constraints, 
the leaves $N$ of $R$ satisfy $\Pi^{\sharp}TN ^{\circ}\oplus TN=TM|_N$,
which by definition means that they are cosymplectic submanifolds.  
There is an induced  Poisson structure on $N$  \cite[Sect. 8]{CrF}, obtained pulling back to $N$ the Dirac structure given by the graph of $\Pi$. The corresponding Poisson bracket of functions $f,g$ on $N$ is 
 $\{\tilde{f},\tilde{g}\}_{\Pi}|_N$, where $\tilde{f},\tilde{g}$ are 
extensions of $f,g$ to $M$ and 
   $d\tilde{f}$ is required to annihilate $\Pi^\sharp(TN ^{\circ})$
at points of $N$.

ii) One checks easily that $\{\varphi^i,g\}_{Dirac}=0$ for all $g\in
C^{\infty}(U)$, i.e. that the $\varphi^i$ are Casimir functions for
$\Pi_{Dirac}$, hence the level sets $N$ of $\varphi$
are Poisson submanifolds (i.e.  unions of symplectic leaves)
w.r.t. $\Pi_{Dirac}$. The Poisson structure on $N$ as a Poisson submanifold of  (M,$\Pi_{Dirac}$) agrees with  the one induced by 
$\Pi$ in the way described in i). Indeed
for all functions $f,g$ on $N$ and extensions $\tilde{f},\tilde{g}$ as above we have 
  $\{\tilde{f},\tilde{g}\}_{Dirac}|_N=\{\tilde{f},\tilde{g}\}_{\Pi}|_N$ since  
$\{\tilde{f}, \varphi^a\}_{\Pi}|_N=0$ for all constraints $\varphi^a$. 
\end{proof}

Now consider an integrable distribution ${R}\subset TM$ and let 
  $D$ be a Dirac structure
 on $E_0\rightarrow M$ so that $D\cap  R ^{\circ}$  has constant rank\footnote{Here   
 $ R ^{\circ}\subset T^*M$ denotes
the annihilator of $R$, i.e., the sections of $ R ^{\circ}$ are the 1-forms
that kill all sections of $R$.}.
 We consider the stretched Dirac subbundle
 $D^{{ R ^{\circ}}}$.
 
The following proposition shows that in the special case that $D$ is the graph of a Poisson structure $\Pi$ and the leaves of $R$ are cosymplectic in $(M,\Pi)$, the Dirac subbundle $D^{{ R ^{\circ}}}$ gives exactly the classical Dirac bracket  (Def. \ref{Dirbr}). Hence $D^{{ R ^{\circ}}}$ can be considered as a generalization of the classical Dirac bracket.

 \begin{proposition}\label{equiv}
1) $D^{ R ^{\circ}}$ is a Dirac structure. It is constructed putting
together the integral submanifolds $N$ of $R$, with the (smooth) Dirac structure induced  pulling back  $D$.

 \vspace{1mm}
\noindent Now assume that  $D$ is the graph of a  Poisson structure $\Pi$.

2a)  $D^{{ R ^{\circ}}}$ 
is itself the graph of   a  Poisson structure iff  the leaves of the distribution $ R$ are  
\emph{Poisson-Dirac submanifolds} \cite{CrF} of $(M,\Pi)$.

2b)
Suppose the stronger condition that the leaves of $R$ are
\emph{cosymplectic} submanifolds of $(M,\Pi)$, 
 so that the Dirac bracket \eqref{dir} can be
defined (see Lemma \ref{PAM}).
 Then $D^{{ R ^{\circ}}}$ is the graph of the  Poisson structure $\Pi_{Dirac}$.
\end{proposition}  
  
\begin{proof}

1) Notice that since $R$ is integrable we can choose a frame for $R^{\circ}$
consisting of closed 1-forms, which act trivially under the Courant bracket.
Hence $R^{\circ}$ is canonical for $D$ (see Def. \ref{canonical}). So from Prop.
 \ref{imply}  we conclude  that $D^{ R ^{\circ}}$ is a Dirac structure.  
 
Since  $\pi(D^{{ R ^{\circ}}})$ is everywhere tangent to the foliation
$R$, the integral submanifolds $N$ of $R$ are unions of presymplectic leaves of 
$D^{{ R ^{\circ}}}$. The Dirac structure $D$ can be restricted to any leaf $N$ of the foliation induced by $R$, delivering a smooth subbundle since $D\cap  R ^{\circ}$ has constant rank \cite{Cou} 
\cite{BW}.
Further, a simple computation shows that 
the pullback to $N$ of    $D$   is equal to
the pullback to $N$ of  $D^{{ R ^{\circ}}}$.   

2a) 
The Dirac structure $D^{{ R ^{\circ}}}$ is the graph of a bivector
field if and only  if
$D^{{ R ^{\circ}}}+TM=TM\oplus T^*M$.
Taking orthogonals  we obtain 
$(D+{ R} ^{\circ})\cap{ R}=\{0\}$.
This can be rewritten as $\Pi^\sharp({ R} ^{\circ})\cap{ R}=\{0\}$, which by definition \cite[Sect. 8]{CrF} means that the leaves of $R$ are Poisson-Dirac submanifolds of $(M,\Pi)$. 


2b) This follows from Lemma \ref{PAM} ii) and part i) of this Proposition.
 \end{proof}
 
 \begin{remark}\label{second}  Prop. \ref{equiv} 2a) shows that
even within the framework of Poisson geometry, i.e. in the case that both $D$ and $D^{ R ^{\circ}}$
 correspond to Poisson structures, our construction of Dirac structure $D^{ R ^{\circ}}$
is more general than the classical
 Dirac bracket.   
   \end{remark}
   
\subsection {The Marsden-Ratiu reduction.} 

We show that the reduced Poisson structure induced via Marsden-Ratiu reduction \cite{MarRat} from a Poisson manifold $(M,\Pi)$  is obtained pushing forward not $\Pi$ itself but rather a suitable stretching of $\Pi$.

We start by recalling the Poisson reduction by distributions as it was stated
by Marsden and Ratiu in \cite{MarRat}, see also \cite{OrtRat}.
The set-up  is the following:
 
\begin{quotation}
$(M,\{\cdot,\cdot\})$ is a Poisson manifold,\\
$N$ is a submanifold with embedding $\iota: N\hookrightarrow M$,\\
${B}\subset T_N M$ is  a smooth
subbundle of $TM$ restricted to $N$.
\end{quotation} 
\noindent We shall also assume that 
${{F}}:={B}\cap TN$ is an integrable regular distribution on $N$
and 
$\un{N}:=N/F$ is a smooth manifold.
 

 
\begin{definition}\label{def:reducible}\cite{MarRat}
 $(M,\{\cdot,\cdot\}, {N}, B)$ is {\bf Poisson reducible}
if there is a Poisson bracket $\{\cdot,\cdot\}_{\un{N}}$ on $\un{N}$
such that for any
$f_1, f_2\in  C^\infty(\un{N})\cong  C^\infty(N)_{{F}}$ we have:
$$\{f_1,f_2\}_{\un{N}}=\iota^*\{f_1^B,f_2^B\}$$
for all extensions $f_i^B\in C^\infty(M)_{B}$ of $f_i$.

 Here  $C^\infty(N)_{F} :=\{f\in C^\infty(N)\mid
\dd f\vert_{F}=0\}$ and $C^\infty(M)_{B}:=\{f\in C^\infty(M)\mid
\dd f\vert_{B}=0\}$. 
\end{definition}



Given  $(M,\{\cdot,\cdot\}, {N}, B)$ 
clearly there is at most\footnote{References \cite{MarRat}, and subsequently \cite{FZ}, formulate  conditions which ensure that $(M,\{\cdot,\cdot\}, {N}, B)$ is  Poisson reducible.} one Poisson bracket $\{\cdot,\cdot\}_{\un{N}}$ on $\un{N}$ satisfying the requirement of Def. \ref{def:reducible}.
The following proposition (which is essentially \cite[Prop. A.2]{FZ}) describes
 the reduced Poisson structure on $\un{N}$ in terms of bivector fields rather than in terms of brackets: it is obtained from $\Pi$ by \emph{stretching}  along $B$, pulling back to $N$ and then pushing forward to $\un{N}$.
 
\begin{proposition}\label{pushpull}
Assume  that   the prescription  of Def. \ref{def:reducible} gives a well-defined
bivector field on $\un{N}$, denote by $L_{\un{N}}$ its graph,
and denote $L_{\Pi}=graph(\Pi)$.
Then the pullback of the almost Dirac structure  $L_{\un{N}}$ under
$p: N \rightarrow \un{N}$ is $\iota^*(L_{\Pi}^B)$.
 
Consequently  $L_{\un{N}}$ is given by the 
  push-forward under $p$ of $\iota^*(L_{\Pi}^B)$.\end{proposition}

\subsection{Couplings on Poisson fibrations} 
Given a Dirac subbundle $D$ and an
isotropic subbundle $S$, there are situations in which one wants to ``deform'' $D$
to a new Dirac subbundle which contains $S$. A natural candidate for the new
Dirac structure is the stretching $D^S$. An instance is provided by our next example,  inspired  by the results of Brahic and Fernandes
\cite{BF} which (in the case of a flat connection) can be rephrased in 
our formalism.

Our data are a manifold $M$ and: 
\begin{quotation}
A splitting of the tangent bundle into two regular, integrable 
distributions: $TM=Hor\oplus Vert$.\\
A two form in $Hor$: $\omega\in\Omega^2(Hor)$.\\
A bivector field in $Vert$: $\pi_V\in\wedge^2(Vert)$.
\end{quotation}
The question is how two combine these data and which are the conditions
that produce a Dirac structure. In principle there are two dual ways
of doing this by using the deformation by stretching. The two different
procedures give the same result.

\vskip 2mm 
\noindent \textbf{a)} Consider $\pi_V^\sharp:T^*M\longrightarrow Vert$ and take
$D=\rm{graph(\pi^\sharp_V)}$.

D is a Dirac structure if and only if

\vskip 1mm
\noindent $i)$ $[\pi_V,\pi_V]=0$.
\vskip 1mm

To define $S$, the stretching direction, consider the bundle map
$\hat\omega^\flat:Hor\longrightarrow Vert^{\circ}$
induced by $\omega$ and take $S={\rm graph}(\hat\omega^\flat)$.

$S$ is involutive if and only if:

\vskip 1mm
\noindent $ii)$ $\omega$ is horizontally closed,
\vskip 1mm
\noindent $iii)$ ${\mathcal L}_v(\omega(u_1,u_2))=0$ for $v\in\Gamma(Vert)$
and $u_i\in\Gamma(Hor)$ s. t. $[v,u_i]\in\Gamma(Vert)$.
\vskip 1mm

Now, assuming that the other conditions hold,  we can show that $S$ 
preserves $D^S$ if and only if

\noindent $iv)$ ${\mathcal L}_u \pi_V=0$ for any $u\in\Gamma(Hor)$ s. t.
$[v,u]\in\Gamma(Vert),\ \forall\  v\in\Gamma(Vert)$.\hfill\break

If conditions $i)$-$iv)$ are satisfied then, using Thm. \ref{connection} ($S^\perp= S + Vert +Hor^\circ$ and therefore $\pi(S^\perp)=TM$), we have that $D^S$ defines a Dirac structure. In the next paragraph
we shall show an alternative way of obtaining the same result with
the roles of $\omega$ and $\pi_V$ exchanged.

\vskip 2mm
\noindent \textbf{b)} We introduce first the bundle map 
$\omega^\flat:Hor\longrightarrow Hor^*$.

Consider the Dirac subbundle $D'\subset TM\oplus T^*M$ induced by
${\rm graph}(\omega^\flat)$, i.e.
$$D'=\{(v,\xi)\vert   v\in Hor, \xi|_{Hor}=\omega^\flat v \}.$$

One can show that $D'$ is a Dirac structure if and only if
condition $ii)$ above holds.

Now we proceed to define the new stretching subbundle $S'$.
Take $$S'_x=\{((\pi_V\xi)_x, \xi_x)\vert \xi_x\in(Hor^{\circ})_x\}.$$

$S'$ is involutive if and only if conditions $i)$ and $iv)$ hold.

Finally one can show that, assuming all previous conditions, 
the stretching $D'^{S'}$ is a 
Dirac structure if and only if 

\vskip 1mm
\noindent $iii)'$ ${\mathcal L}_v(\omega(u_1,u_2))=0$ for 
$v\in\Gamma(\pi_V^\sharp(T^*M))$ and $u_i\in\Gamma(Hor)$ s. t. 
$[v,u_i]\in\Gamma(Vert)$.
\vskip 1mm

It is interesting to compare the two construction.
First it is clear that $D^S=D'^{S'}$. 
Further, condition $iii)$ in construction a)
implies condition $iii)'$ in b). Therefore, even if both give the same
final result, the second construction has a broader 
range of application.

\begin{remark} 
We establish the connection between the above and  the coupling of Poisson
fibrations of ref. \cite{BF}.
Suppose that $M$ is the total space of a fibration so that $Vert$ is the distribution tangent to the fibers.
One computes easily that $D^S$ agrees with the fiber non-degenerate  
almost Dirac structure associated to the triple $(\pi_V,Hor,\omega)$ 
in Cor. 2.6 of \cite{BF}.
Brahic and Fernandes compute the necessary and sufficient conditions for 
this to be a Dirac structure in Cor. 2.8 of \cite{BF}.
If the horizontal connection is flat their conditions
are equivalent to our $i), ii), iii)'$ and $iv)$ above,
i.e. the conditions for having a Dirac structure following
the stretching procedure introduced in the paper.
\end{remark}

\vspace{5mm} 
  
\noindent{\bf Acknowledgments:} Research partially supported by
grants FPA2003-02948 and FPA2006-02315, MEC (Spain) and by the
Forschungskredit 2006 of the University of Z\"urich.
 

\providecommand{\bysame}{\leavevmode\hbox to3em{\hrulefill}\thinspace}
\providecommand{\MR}{\relax\ifhmode\unskip\space\fi MR }
\providecommand{\MRhref}[2]{%
  \href{http://www.ams.org/mathscinet-getitem?mr=#1}{#2}
}
\providecommand{\href}[2]{#2}

\end{document}